\documentclass[11pt,a4paper]{article}
\usepackage{epsfig,amsfonts,amsgen,amsmath,amstext,amsbsy,amsopn,amsthm}
\usepackage{color}
\newtheorem{false statement}{False statement}

\theoremstyle{definition}

\newcounter{mathitem}
  {\begin{list}{{$(\roman{mathitem})$}}{
   \setcounter{mathitem}{0}
   \usecounter{mathitem}
   \setlength{\topsep}{0pt plus 2pt minus 0pt}
   \setlength{\parskip}{0pt plus 2pt minus 0pt}
   \setlength{\partopsep}{0pt plus 2pt minus 0pt}
   \setlength{\parsep}{0pt plus 2pt minus 0pt}
   \setlength{\leftmargin}{35pt}
   \setlength{\itemsep}{0pt plus 2pt minus 0pt}}}
  {\end{list}}

\begin{document}

\title
{\bf A note on nowhere-zero 3-flow and $Z_3$-connectivity}

\author{
Fuyuan Chen$^a$\thanks{E-mail address: chenfuyuan19871010@163.com (F. Chen).}
and Bo Ning$^{b}$\thanks{E-mail address: ningbo\_math84@mail.nwpu.edu.cn (B. Ning).}\\
\small $^{a}$Center for Discrete Mathematics,\\
\small Fuzhou University, Fuzhou, Fujian 350108, P.R.~China\\
\small $^{b}$Department of Applied Mathematics, School of Science,\\
\small  Northwestern Polytechnical University, Xi'an, Shaanxi 710072,
P.R.~China\\[2mm]}
\maketitle

\noindent\bf ABSTRACT. \rm  There are many major open problems in integer flow theory, such as Tutte's 3-flow
conjecture that every 4-edge-connected graph admits a nowhere-zero 3-flow, Jaeger et al.'s conjecture that every 5-edge-connected graph is
$Z_3$-connected and Kochol's conjecture that every bridgeless graph with at most three 3-edge-cuts admits a nowhere-zero 3-flow (an equivalent version of 3-flow conjecture). Thomassen proved that every 8-edge-connected graph is $Z_3$-connected and therefore admits a nowhere-zero 3-flow. Furthermore, Lov$\acute{a}$sz, Thomassen, Wu and Zhang improved Thomassen's result to 6-edge-connected graphs. In this paper, we prove that: (1) Every 4-edge-connected graph with at most seven 5-edge-cuts admits a nowhere-zero 3-flow. (2) Every bridgeless graph containing no 5-edge-cuts but at most three 3-edge-cuts admits a nowhere-zero 3-flow. (3) Every 5-edge-connected graph with at most five 5-edge-cuts is $Z_3$-connected. Our main theorems are partial results to Tutte's 3-flow conjecture, Kochol's conjecture and Jaeger et al.'s conjecture, respectively.
\section{Introduction}
All graphs considered in this paper are loopless, but allowed to have multiple edges. A graph $G$ is called $k$-edge-connected, if $G - S$ is connected for each edge set $S$ with $|S|< k$. Let $X$, $Y$ be two disjoint subsets of $V(G)$. Let $\partial_{G}(X, Y)$ be the set of edges of $G$ with one end in $X$ and the other in $Y$. In particular, if $Y= \overline{X}$, we simply write $\partial_{G}(X)$ for $\partial_{G}(X, Y)$, which is the {\em edge-cut} of $G$ associated with $X$. The edge set $C=\partial_{G}(X)$ is called a {\em $k$-edge-cut} if $|\partial_{G}(X)|=k$. If $X$ is nontrivial, we use $G/ X$ to denote the graph obtained from $G$ by replacing $X$ by a single vertex $x$ that is incident with all the edges in $\partial_{G}(X)$.

Let $D$ be an orientation of $E(G)$. The {\em out-cut} of $D$ associated with $X$, denoted by $\partial_{D}^{+}(X)$, is the set of arcs of $D$ whose tails lie in $X$. Analogously, the {\it in-cut} of $D$ associated with $X$, denoted by $\partial_{D}^{-}(X)$, is the set of arcs of $D$ whose heads lie in $X$. We refer to $|\partial_{D}^{+}(X) |$ and $|\partial_{D}^{-}(X)|$ as the out-degree and in-degree of $X$, and denote these quantities by $d_{D}^{+}(X)$ and $d_{D}^{-}(X)$, respectively.

\vskip10pt\noindent{\bf Definition 1.1.} \sl
(1) An orientation $D$ of $E(G)$ is called a {\em modulo 3-orientation} if $$d_{D}^{+}(v) - d_{D}^{-}(v) \equiv 0 \pmod 3$$ for every vertex $v\in V(G)$.\\
(2) A pair $(D, f)$ is called a {\em nowhere-zero 3-flow} of $G$ if $D$ is an orientation of $E(G)$ and $f$ is a function from $E(G)$ to $\{\pm 1, \pm2\}$, such that $$\sum_{e\in \partial_{D}^{+}(v)}f(e)=\sum_{e\in \partial_{D}^{-}(v)}f(e)$$ for every vertex $v\in V(G)$.\rm

The 3-flow conjecture, proposed by Tutte as a dual version of Gr\"{o}tzsch's 3-color theorem for planar graphs, may be one of the most major open problems in integer flow theory.

\vskip10pt\noindent {\bf Conjecture 1.2.} \sl (3-Flow conjecture, Tutte \cite{Tu}) Every 4-edge-connected graph admits a nowhere-zero
3-flow.\rm \vskip10pt

Kochol proved that Tutte's 3-flow conjecture is equivalent to the following two conjectures.

\vskip10pt\noindent {\bf Conjecture 1.3.} (Kochol \cite{K}) \sl Every
$5$-edge-connected graph admits a nowhere-zero 3-flow.\rm\vskip10pt

\vskip10pt\noindent {\bf Conjecture 1.4.} (Kochol \cite{K1}) \sl Every bridgeless graph with at most three 3-edge-cuts admits a nowhere-zero 3-flow.\rm\vskip10pt

A weakened version of Conjecture 1.2, the so-called weak 3-flow conjecture, was proposed by Jaeger.

\vskip10pt\noindent {\bf Conjecture 1.5.} (Weak 3-flow conjecture, Jaeger \cite{J}) \sl There is a natural number $h$ such that every $h$-edge-connected graph admits a nowhere-zero 3-flow.\rm\vskip10pt

Lai and Zhang \cite{LZ} and Alon et al. \cite{ALM} gave partial results on Conjectures 1.2 and 1.5.

\vskip10pt\noindent {\bf Theorem 1.6.} (Lai and Zhang \cite{LZ}) \sl Every $4\lceil \log_{2}n_0\rceil$-edge-connected graph with at most $n_0$ odd-degree vertices admits a nowhere-zero 3-flow.\rm\vskip10pt

\vskip10pt\noindent {\bf Theorem 1.7.} (Alon, Linial and Meshulam \cite{ALM}) \sl Every $2\lceil \log_{2}n\rceil$-edge-connected graph with $n$ vertices admits a nowhere-zero 3-flow.\rm\vskip10pt

Recently, Thomassen \cite{T} confirmed weak 3-flow conjecture. He proved

\vskip10pt\noindent {\bf Theorem 1.8.} (Thomassen \cite{T}) \sl Every 8-edge-connected graph is $Z_3$-connected and therefore admits a nowhere-zero 3-flow.\rm\vskip10pt

Thomassen's method was further refined by Lov\'{a}sz, Thomassen, Wu and Zhang \cite{LTWZ} to obtain the following theorem.

\vskip10pt\noindent {\bf Theorem 1.9.} (Lov\'{a}sz, Thomassen, Wu and Zhang \cite{LTWZ}) \sl Every 6-edge-connected graph is $Z_3$-connected and therefore admits a nowhere-zero 3-flow.\rm\vskip10pt

For more results on Tutte's 3-flow conjecture, we refer the reader to the introduction part of \cite{LTWZ} and the book written by Zhang \cite{Z}.

In this paper, we will give the following conjecture which is equivalent to Tutte's 3-flow conjecture.

\vskip10pt\noindent {\bf Conjecture 1.10.} \sl Every 5-edge-connected graph with minimum degree at least $6$ has a nowhere-zero 3-flow.\rm

To prove the equivalence of Conjectures 1.2 and 1.10, the following lemma is needed.

\vskip10pt\noindent {\bf Lemma 1.11.} (Tutte \cite{Tu1}) \sl  Let {\em $F(G, k)$} be the number of
nowhere-zero $k$-flows of $G$. Then $F(G, k)= F(G /e, k)- F(G\setminus e, k)$ if $e$ is not a loop of $G$.\rm

\vskip10pt\noindent{\bf Proposition}. \sl Conjectures 1.2 and 1.10 are equivalent.\rm
\vskip5pt\noindent{\bf Proof.} It is obvious that
Conjecture 1.2 implies Conjecture 1.3, and Conjecture 1.3 implies Conjecture 1.10. Now we prove that
Conjecture 1.10 can imply Conjecture 1.3. Let $G$ be a
5-edge-connected graph. Let $G'$ be the graph obtained from $G$ by gluing $|V(G)|$ disjoint copies of $K_7$, such that for each such copy $H_i$, $|V(H_i)\cap V(G)|=1$ ($i=1, 2, \cdots, |V(G)|$). Then $G'$ is
5-edge-connected and its minimum degree is at least $6$, and thus has a nowhere-zero 3-flow. By Lemma 1.11, $G$ has a nowhere-zero 3-flow. Therefore
Conjecture 1.10 implies Conjecture 1.3. Note that Conjecture 1.2 is equivalent to Conjecture 1.3. This completes the proof. \quad\vrule height6pt
width6pt depth0pt

Our first main result is the following theorem.

\vskip10pt\noindent{\bf Theorem 1.12.} \sl Let $G$ be a bridgeless graph and let $P=\{C=\partial_{G}(X): |C|=3,\ X\subset V(G)\}$ and $Q=\{C=\partial_{G}(X): |C|=5,\  X\subset V(G)\}$. If $2|P|+|Q|\leq 7$, then $G$ has a modulo 3-orientation (and therefore
has a nowhere-zero 3-flow).\rm

As corollaries of Theorem 1.12, we obtain Theorems 1.13 and 1.14.

\vskip10pt\noindent{\bf Theorem 1.13.} \sl Every 4-edge-connected graph with at most seven 5-edge-cuts admits a nowhere-zero 3-flow.\rm

\vskip10pt\noindent{\bf Theorem 1.14.} \sl Every bridgeless graph containing no 5-edge-cuts but at most three 3-edge-cuts admits a nowhere-zero 3-flow.\rm

\vskip10pt\noindent{\bf Remark.} \sl
The number of 3-edge-cuts in Theorem 1.14 can not be improved from three to four since $K_4$ or any graph contractible to $K_4$ has no nowhere-zero 3-flow.\rm

Theorems 1.13 and 1.14 partially confirm Conjectures 1.2 and 1.4, respectively.

\vskip10pt\noindent{\bf Definition 1.15.} \sl
(1) A mapping $\beta_{G}: V(G) \mapsto Z_k$ is called a $Z_k$-boundary of $G$ if $$\sum_{v \in V(G)} \beta_{G}(v) \equiv 0
 \pmod{k}$$

(2) A graph $G$ is called {\em $Z_k$-connected}, if for every $Z_k$-boundary $\beta_{G}$, there is an orientation $D_{\beta_{G}}$ and a function $f_{\beta_{G}}$:
$E(G)\mapsto Z_k - \{0\}$, such that $$\sum_{e\in \partial_{D_{\beta_{G}}}^{+}(v)} f_{\beta_{G}}(e) - \sum_{e\in \partial_{D_{\beta_{G}}}^{-}(v)} f_{\beta_{G}}(e) \equiv \beta_{G}(v) \pmod k$$ for every vertex $v\in V(G)$. \rm

Jaeger, Linial, Payan and Tarsi \cite{JLPT} conjectured that

\vskip10pt\noindent {\bf Conjecture 1.16.} (Jaeger, Linial, Payan and Tarsi \cite{JLPT}) \sl Every
$5$-edge-connected graph is $Z_3$-connected.\rm\vskip10pt

By applying a similar argument in the proof of Theorem 1.12, we could obtain the second main result which is a partial result to Conjecture 1.16.

\vskip10pt\noindent{\bf Theorem 1.17.} \sl Every 5-edge-connected graph with at most five 5-edge-cuts is $Z_3$-connected.\rm

In the next section, some necessary preliminaries will be given. In Sections 3 and 4, proofs of Theorems 1.12 and 1.17 will be given, respectively.
\section{Preliminaries}
\label{S: Preliminaries}
In this section, we will give additional but necessary notations and definitions, and then give some useful lemmas.

\vskip10pt\noindent{\bf Definition 2.1} \sl Let $\beta_{G}$ be a $Z_3$-boundary of $G$. An orientation $D$ of
$G$ is called a {\em $\beta_{G}$-orientation} if $$d^+_D(v) - d^-_D(v) \equiv \beta_{G}(v) \pmod{3}$$ for every vertex $v
\in V(G)$.\rm

Let $G$ be a graph and $A$ be a vertex subset of $G$. The {\em
degree} of $A$, denoted by $d_{G}(A)$, is the number of edges with
precisely one end in $A$. Moreover if $A=\{x\}$, we simply write
$d_{G}(x)$.

Let $G$ be a graph and $\beta_{G}$ be a $Z_3$-boundary of $G$. Define a
mapping $\tau_{G}:V(G)\mapsto \{ 0, \pm 1, \pm 2, \pm 3\}$ such that,
for each vertex $x\in V(G)$,
\begin{equation}
\tau_{G}(x) \equiv \left\{ \begin{array}{rl}
\beta_{G} (x) & \pmod 3 \\
 d_{G}(x) & \pmod 2.
 \end{array} \right.\nonumber
\label{EQ: d t}
\end{equation}
Now, the mapping $\tau_{G}$ can be further extended to any nonempty
vertex subset $A$ as follows: $$ \tau_{G}(A) \equiv \left\{
\begin{array}{rl}
\beta_{G} (A) & \pmod 3 \\
 d_{G}(A) & \pmod 2.
 \end{array} \right.
$$
where $\beta_{G} (A)\equiv\sum_{x\in A}\beta_{G} (x) \in\{0,1,2\} \pmod{3}$.

\vskip10pt\noindent{\bf Proposition 2.2} \sl Let $G$ be a graph and $A$ be a vertex subset
of $G$.

(1) If $d_{G} (A)\leq 5$, then $d_{G} (A)\leq 4+|\tau_{G} (A)|$.

(2) If $d_{G} (A)\geq 6$, then $d_{G} (A)\geq 4+|\tau_{G} (A)|$.\rm

Proposition 2.2 follows from the fact that $|\tau_{G}(A)|\leq 3$ and
$d_{G}(A)-|\tau_{G}(A)|$ is even.\rm

\vskip10pt\noindent{\bf Lemma 2.3 (Tutte \cite{Tu})} \sl Let $G$ be a graph.

(1) $G$ has a nowhere-zero 3-flow if and only if $G$
has a modulo 3-orientation.

(2) $G$ has a nowhere-zero 3-flow if and only if $G$ has a $\beta_{G}$-orientation with
$\beta_{G}=0$.\rm

The following lemma is Theorem 3.1 in \cite{LTWZ} by Lov\'asz et al. This lemma will play the main role in our proofs.

\vskip10pt\noindent{\bf Lemma 2.4 (Lov\'asz, Thomassen, Wu and Zhang \cite{LTWZ})}\sl Let $G$ be a graph, $\beta_{G}$ be a $Z_3$-boundary of $G$, and let
$z_0\in V(G)$ and $D_{z_0}$ be a pre-orientation of $E(z_0)$ of all edges incident with $z_0$. Assume
that

(i) $|V(G)|\geq 3$.

(ii) $d_{G} (z_0)\leq 4+|\tau_{G}(z_0)|$ and $d_{D_{z_0}}^+ (z_0)-d_{D_{z_0}}^- (z_0)\equiv \beta_{G} (z_0)\pmod
 3$, and

(iii) $d_{G}(A)\geq 4+|\tau_{G} (A)|$ for each
nonempty vertex subset $A$ not containing $z_0$ with $|V
(G)\setminus A|>1$.

Then the pre-orientation $D_{z_0}$ of $E (z_0)$ can be extended to
an orientation $D$ of the entire graph $G$, that is, for every
vertex $x$ of $G$,
$$d_{D}^+ (x)-d_{D}^- (x)\equiv \beta_{G} (x) \pmod 3.$$\rm

\section{Proof of Theorem 1.12}

If not, suppose that $G$ is a counterexample, such that $|V(G)|+|E(G)|$ is as small as possible.
Let $P'=\{x\in V(G): d_{G}(x)=3\}$ and $Q'=\{x\in V(G): d_{G}(x)=5\}$.

\noindent{\bf Claim 1.} $|V(G)|\geq 3$.

\vskip3pt\noindent{\bf Proof}.\noindent \ If $|V(G)|= 1$, then $G$ has a nowhere-zero 3-flow, a contradiction. If $|V(G)|=2$, let
$V(G)=\{x,y\}$, then all the edges of $G$ are all between $x$ and $y$. Since $G$ is bridgeless, $|E(G)|\geq 2$. Let $a$ be the integer in $\{0,1,2\}$ such that $a\equiv |E(G)|-a \pmod 3$. Orient $a$ edges from $x$ to $y$ and the remaining $|E(G)|-a$ edges from $y$ to $x$. Clearly, the resulting orientation is a modulo $3$-orientation of $G$, a contradiction. Therefore $|V(G)|\geq 3$.\quad\vrule height6pt
width6pt depth0pt

\noindent{\bf Claim 2.} $G$ is 3-edge-connected, and $G$ has no nontrivial 3-edge-cuts.

\vskip3pt\noindent{\bf Proof}.\noindent\ If $G$ has a vertex $x$ of degree 2, then suppose that $xx_1, xx_2\in E(G)$. By the minimality of $G$, $(G-\{xx_1, xx_2\})\cup \{x_1x_2\}$ has a nowhere-zero 3-flow $f'$. However, $f'$ can be extended to a nowhere-zero 3-flow $f$ of $G$, a contradiction. If $G$ has a nontrivial $k$-edge-cut$(k=2, 3)$, then contract one side and find a mod 3-orientation by the minimality of $G$. Merge such two mod 3-orientations and we will get one for $G$, a contradiction.\quad\vrule height6pt
width6pt depth0pt

\noindent{\bf Claim 3.} For any $U\subset V(G)$, if $d_{G}(U)\leq 5$ and $|U|\geq 2$, then $U\cap (P'\cup Q')\neq\emptyset$.

\vskip3pt\noindent{\bf Proof}.\noindent\ If not, choose $U$ to be a minimal one such that: for any $A\subset U$ with
$2\leq |A|<|U|$, we have $d_{G}(A)\geq 6$.

By the minimality of $G$, $G/U$ has a modulo 3-orientation $D'$ which is a partial modulo
3-orientation of $G$, such that $d_{D'}^+(x)\equiv d_{D'}^-(x)\pmod 3$ for
each $x\in V(G)\setminus U$.

Let $G'$ be a graph obtained from $G$ by contracting $V(G)\setminus U$ as $z_0$ and let $\beta_{G'}=0$.

(i) Since $V(G')=U+z_0$, $|V(G')|=|U|+1\geq 3$.

(ii) Since $d_{G'}(z_0)=d_{G}(U)\leq 5$, by Proposition 2.2 (1), $d_{G'}(z_0)\leq 4+|\tau_{G'}(z_0)|$.

(iii) By the assumption and minimality of $U$, we have that $\forall A\subset U$, $d_{G}(A)\neq 5$ and $d_{G}(A)\neq 3$. If $d_{G}(A)=4$, then $d_{G'}(A)=d_{G}(A)=4$ and $\tau_{G'}(A)=\beta_{G'}(A)=\beta_{G}(A)=0$. Thus $d_{G'}(A)=4=4+|\tau_{G'}(A)|$. If $d_{G}(A)\geq 6$, then by Proposition 2.2 (2), $d_{G'}(A)=d_{G}(A)\geq 4+|\tau_{G'}(A)|$.

By Lemma 2.4, we could see that the pre-orientation of $E'(z_0)$ of all edges incident with $z_0$ can be extended to a $\beta_{G'}$-orientation of $G'$. Then $G$ has a modulo 3-orientation, which is a contradiction. \quad\vrule height6pt
width6pt depth0pt

Let $G_1'$ be a graph obtained from $G$ by adding a new vertex $z_0$ and $2|P'|+|Q'|$ edges between $z_0$ and $P'\cup Q'$, such that:

(i) For each vertex $v\in P'$, we add two arcs with the same direction between it and $z_0$; and

(ii) For each vertex $v\in Q'$, we add one arc between it and $z_0$.

If $2|P'|+|Q'|\leq 5$, then all added arcs could be from $z_0$ to $P'\cup Q'$. Define $\beta_{G_1'}$ as follows:

(1) $\beta_{G_1'}(x)=0$ if $x\not\in (P'\cup Q')+z_0$;

(2) $\beta_{G_1'}(x)=1$ if $x\in P'$;

(3) $\beta_{G_1'}(x)=2$ if $x\in Q'$;

(4) $\beta_{G_1'}(z_0)\equiv 2|P'|+|Q'|\pmod 3$ and $\beta_{G_1'}(z_0)\in \{0,1,2\}$.

If $2|P'|+|Q'|=6$ or 7, in this case, if $|P'|\neq 0$, choose one vertex $v\in P'$, such that the two arcs with ends $z_0$ and $v$ are from $v$ to $z_0$, the other arcs incident with $z_0$ are all directed from $z_0$. If $|P'|=0$, then two arcs are from $Q'$ to $z_0$, the others verse. Define $\beta_{G_1'}$ as follows:

(1) $\beta_{G_1'}(x)=0$ if $x\not\in (P'\cup Q')+z_0$;

(2) $\beta_{G_1'}(x)=2$ if $x\in Q'$ and the arc $(z_0, x)$ exists or $x\in P'$ and the two arcs with ends $z_0$ and $x$ are from $x$ to $z_0$;

(3) $\beta_{G_1'}(x)=1$ if $x\in Q'$ and the arc $(x, z_0)$ exists or $x\in P'$ and the two arcs with ends $z_0$ and $x$ are from $z_0$ to $x$;

(4) $\beta_{G_1'}(z_0)\equiv (2|P'|+|Q'|-2)-2\pmod 3$.

Now $d_{G_1'}(z_0)\leq 4+|\tau_{G_1'}(z_0)|$ and $|V(G_1')|=|V(G)|+1\geq 4$. We claim that: $d_{G_1'}(A)\geq 4+|\tau_{G_1'}(A)|$, for each
nonempty vertex subset $A$ not containing $z_0$ with $|V(G_1')\setminus A|>1$.

If $A\cap (P'\cup Q')=\emptyset$, then by Claim 3, $d_{G}(A)=4$ or $d_{G}(A)\geq 6$. In each case we could get that $d_{G_1'}(A)=d_{G}(A)\geq 4+|\tau_{G_1'}(A)|$.

If $A\cap (P'\cup Q')\neq\emptyset$, then by Claim 2, $d_{G_1'}(A)\geq 5$. If $d_{G_1'}(A)=5$, then $d_{G}(A)=3$ or 4 and $|A\cap (P'\cup Q')|=1$, and it follows that $\beta_{G_1'}(A)=1$ or 2, and $|\tau_{G_1'}(A)|=1$. Then $d_{G_1'}(A)\geq 4+|\tau_{G_1'}(A)|$. If $d_{G_1'}(A)\geq 6$, by Proposition 2.2 (2), we have that $d_{G_1'}(A)\geq 4+|\tau_{G_1'}(A)|$.

Now $G_1'$ satisfies all the conditions of Lemma 2.4. By Lemma 2.4, $G_1'$ has a $\beta_{G_1'}$-orientation extended from the pre-orientation of $E_1'(z_0)$ of all edges incident with $z_0$, which implies that $G$ has a $\beta_{G}$-orientation with $\beta_{G}=0$. By Proposition 2.3, $G$ has a nowhere-zero 3-flow, which is a contradiction.{\hfill$\Box$}

\section{Proof of Theorem 1.17}

Assume not. Suppose that $G$ is a counterexample, such that $|V(G)|+|E(G)|$ is as small as possible. Let $S'=\{x\in V(G): d_{G}(x)=5\}$ and $S=\{C=\partial_{G}(X): |C|=5,\ X\subset V(G)\}$. Let $\beta_{G}$ be a $Z_3$-boundary, such that $G$ has no $\beta_{G}$-orientation.

\noindent {\bf Claim 1.} $|V(G)|\geq 3$ and $|S'|\leq |S|\leq 5$.

\vskip3pt\noindent{\bf Proof}.\noindent\ Since $G$ is 5-edge-connected, $|V(G)|\geq 2$. If $|V(G)|=2$, let
$V(G)=\{x,y\}$. Then all the edges of $G$ are between $x$ and $y$,
and $|E(G)|\geq 5$. Let $D_x$ be an orientation of $x$, such that $d_{D_x}^+(x)-d_{D_x}^-(x)\equiv \beta_{G}(x) \pmod 3$. Since $\beta_{G}$ is a $Z_3$-boundary, $d_{D_x}^+(y)-d_{D_x}^-(y)\equiv \beta_{G}(y) \pmod 3$. Therefore $G$ has a $\beta_{G}$-orientation, a contradiction.
Hence $|V(G)|\geq 3$ and $|S'|\leq |S|\leq 5$.\quad\vrule height6pt
width6pt depth0pt

\noindent {\bf Claim 2.} Let $U\subset V(G)$ with  $|U|\geq 2$. If
$d_{G}(U)=5$, then $U\cap S'\neq\emptyset$.

\vskip3pt\noindent{\bf Proof}.\noindent\ If not, choose $U$ to be a minimal one such that: for any $A\subset U$ with
$2\leq |A|<|U|$, we have $d_{G}(A)\neq 5$.

By the minimality of $G$, $G/U$ has a $\beta_{G}$-orientation $D'$ which is a partial $\beta_{G}$-orientation of $G$, such that $d_{D'}^+(x) - d_{D'}^-(x) \equiv \beta_{G}(x) \pmod 3$ for
each $x\in V(G)\setminus U$.

Let $G'$ be a graph obtained from $G$ by contracting $V(G)\setminus U$ as $z_0$, and let $\beta_{G'}=\beta_{G}$.

(i) Since $V(G')=U+z_0$, $|V(G')|=|U|+1\geq 3$.

(ii) Since $d_{G'}(z_0)=d_{G}(U)=5$, by Proposition 2.2 (1), we have that $d_{G'}(z_0)\leq 4+|\tau_{G'}(z_0)|$.

(iii) By the assumption and minimality of $U$, we have that $\forall A\subset U$, $d_{G}(A)\neq 5$.

Therefore $d_{G}(A)\geq 6$. By Proposition 2.2 (2), $d_{G'}(A)=d_{G}(A)\geq 4+|\tau_{G'}(A)|$.

By Lemma 2.4, the pre-orientation of $E'(z_0)$ of all edges incident with $z_0$ can be extended to a $\beta_{G'}$-orientation of $G'$. Then $G$ has a $\beta_{G}$-orientation, which is a contradiction. \quad\vrule height6pt
width6pt depth0pt

Let $G_1'$ be a graph obtained from $G$ by adding a new vertex $z_0$ and $|S'|$ arcs from $z_0$ to $S'$, such that each vertex in $S'$ has degree 6 in $G_1'$.

Define $\beta_{G_1'}$ as follows:

(1) $\beta_{G_1'}(x)=\beta_{G}(x)$ if $x\not\in S'+z_0$;

(2) $\beta_{G_1'}(x)\equiv \beta_{G}(x)-1 \pmod 3$ if $x\in S'$;

(3) $\beta_{G_1'}(z_0)\equiv |S'|\pmod 3$ and $\beta_{G_1'}(z_0)\in \{0,1,2\}$.

Now $d_{G_1'}(z_0)\leq 4+|\tau_{G_1'}(z_0)|$ and $|V(G_1')|=|V(G)|+1\geq 4$. We claim that $d_{G_1'}(A)\geq 4+|\tau_{G_1'}(A)|$, for each
nonempty vertex subset $A$ not containing $z_0$ with $|V(G_1')\setminus A|>1$.

If $A\cap S'=\emptyset$, then by Claim 2, $d_{G_1'}(A)=d_{G}(A)\neq 5$. Thus $d_{G_1'}(A)\geq 6$. By Proposition 2.2 (2), $d_{G_1'}(A)\geq 4+|\tau_{G_1'}(A)|$.

If $A\cap S'\neq\emptyset$, then $d_{G_1'}(A)\geq d_{G}(A)+1\geq 6$. By Proposition 2.2 (2), we have that $d_{G_1'}(A)\geq 4+|\tau_{G_1'}(A)|$.

Now $G_1'$ satisfies all the conditions of Lemma 2.4. By Lemma 2.4, $G_1'$ has a $\beta_{G_1'}$-orientation extended from the pre-orientation of $E_1'(z_0)$ of all edges incident with $z_0$, which implies that $G$ has a $\beta_{G}$-orientation, a contradiction.

The proof is complete.{\hfill$\Box$}

\section*{Acknowledgement}
This work is supported by NSFC (No.~11271300) and the Doctorate
Foundation of Northwestern Polytechnical University (cx201326).

\end{document}